\documentclass[twocolumn]{autart}  
\usepackage{graphicx}          
\usepackage{amsmath}
\usepackage{amssymb}
\usepackage{bm}
\usepackage{bbm}
\usepackage[sort&compress]{natbib}
\usepackage{xcolor}

\makeatletter
\newcommand{\onetagright}{\tagsleft@false}
\makeatother
\begin{document}

\begin{frontmatter}
\title{Consensus Control for Linear Systems \\with Optimal Energy Cost\thanksref{footnoteinfo}} 

\thanks[footnoteinfo]{This work is supported by China Scholarship Council.}

\author[KTH]{Han Zhang}\ead{hanzhang@kth.se},    
\author[KTH]{Xiaoming Hu}\ead{hu@kth.se},              

\address[KTH]{Department of Mathematics, KTH Royal Institute of Technology, SE-100 44, Stockholm, Sweden}


\begin{abstract}
In this paper, we design an optimal energy cost controller for linear systems asymptotic consensus given the topology of the graph. The controller depends only on relative information of the agents. Since finding the control gain for such controller is hard, we focus on finding an optimal controller among a classical family of controllers which is based on Algebraic Riccati Equation (ARE) and guarantees asymptotic consensus. Through analysis, we find that the energy cost is bounded by an interval and hence we minimize the upper bound. In order to do that, there are two classes of variables that need to be optimized: the control gain and the edge weights of the graph and are hence designed from two perspectives. A suboptimal control gain is obtained by choosing $Q=0$ in the ARE. Negative edge weights are allowed, and the problem is formulated as a Semi-definite Programming (SDP) problem. Having negative edge weights means that "competitions" between the agents are allowed. The motivation behind this setting is to have a better system performance. We provide a different proof compared to \cite{thunberg2016optimal} from the angle of optimization and show that the lowest control energy cost is reached when the graph is complete and with equal edge weights.  Furthermore, two sufficient conditions for the existence of negative optimal edge weights realization are given. 
In addition, we provide a distributed way of solving the SDP problem when the graph topology is regular.
\end{abstract}

\end{frontmatter}

\section{Introduction}
Cooperative control for multi-agent systems has been intensively studied in the past decade, with wide applications on smart power grid, robots and so on. The cooperative control methods usually imply a distributed manner, which make the systems enjoy many advantages such as being robust and economic.

Consensus is an important topic in the research of cooperative control for multi-agent. The goal is to let the states or the outputs of all agents become the same by control laws that depend on the information of the agent and its neighbours. In this paper, we consider the case of asymptotic consensus for linear systems. The goal of this paper is to design distributed controllers using only the relative information between the agents and with minimal control energy cost such that all the system states will eventually become the same as time goes infinity.

In this paper, we focus on designing the optimal energy cost controller for the linear systems so that they can reach asymptotic consensus. It is well-known that the asymptotic consensus for linear systems is equivalent to regulating $N-1$ systems with the dynamics $\dot{x}_i=Ax_i+\lambda_iBu_i$,  $2\leq i \leq N$,  where $\lambda_i$ is the \textit{i}th smallest eigenvalue of the Laplacian matrix, see \citep{fax2004information}, \citep{zhang2011optimal}. The area is well-studied, for example, in \citep{borrelli2008distributed}, the authors design linear quadratic regulators for identical linear systems when the graph topology is given. Though it is a quite similar to what we undertake here, but as mentioned above, the equivalent problem is to regulate $N-1$ different systems with minimum energy, and hence is not the same. The authors of \citep{rogge2010consensus} consider the problem as a quadratic optimal control problem on a ring network while we consider the case of a graph with arbitrary topologies. Augmented Lagrangian approach is used in \citep{lin2011augmented} to design a  structured distributed controller so that the $H_2$ norm of the noisy systems is minimized. \cite{deshpande2011distributed} consider a similar problem and use a two-step approach to design the control law. But their controller does not only use relative information of the agents, but also uses the agents' own states. \cite{cao2010optimal} study the optimal consensus of the single-integrators for both discrete-time and continuous-time case. In \citep{lin2013design}, alternating direction method of multipliers (ADMM) is used to minimize the $H_2$ norm so that an optimal sparse feedback gain is obtained. In \citep{thunberg2016optimal}, a "topology free" control energy minimization problem is considered and the distributed energy-optimal control corresponds to a complete graph with equal edge weights.

On the other hand, distributed optimization has attracted great attentions these years due to the wide applications in the network. Compared to the abundant results on distributed optimization in real vector spaces, for instances, \citep{yi2015distributed}, \citep{annergren2014distributed}, \citep{lou2014approximate}, \citep{nedic2009distributed}, \citep{nedic2010constrained}, the results on distributed SDP still remain limited. \cite{dall2013distributed} considered an optimal power flow problem in power grids and relax it into an SDP problem. Then ADMM is used to solve the problem in a distributed manner. In \citep{simonetto2014distributed}, a sensor localization problem is considered. Relaxation towards an SDP and ADMM is used to solve the problem distributedly as well. \cite{pakazad2014distributed} analyse the robustness of interconnected uncertain systems and linear matrix inequalities (LMIs) are reformulated in SDP. Chordal sparsity structure is assumed among the data matrices of the SDP and hence the problem is decomposed and solved distributedly using proximal splitting method. In \citep{pakazad2015distributed}, coupled SDPs with tree structures are considered. A distributed primal-dual interior point method is proposed to solve the coupled SDPs. The aforementioned work all utilize the idea of "decomposition" somehow but in this paper we treat the SDP in a different manner: reaching an optimal consensus in the intersection of convex feasible sets. Also, what makes our work different from existing distributed optimization problem is that our problem motivates from optimizing a parameter of a graph and the communication network of the distributed optimization algorithm is actually the physical network itself, while most distributed optimization algorithms relax and decompose the original problem and "design" the communication network according to the structure of the decomposed problem, see \citep{pakazad2015distributed} as an example.
 
The main contribution of this paper is the construction of an optimal energy controller that depends only on the relative information between the agents. The controller has two classes of variables that need to be determined: the control gain and the edge weights of the graph. Similar to \citep{borrelli2008distributed}, computing the optimal control gain for the controller is hard, thus we focus on finding an optimal controller among a classical family of controller designs based on ARE and guarantees asymptotic consensus. Through analysis, we found that the energy cost is bounded by an interval and hence we minimize the upper bound. A suboptimal control gain is obtained by choosing $Q=0$ in the ARE; the edge weights of the graph is optimized by solving an underlying SDP. The controller that we designed enjoy several favourable properties:
\begin{enumerate}
\item The controller coincides with the optimal control in \citep{thunberg2016optimal} when the graph is complete. It has been pointed out in \citep{thunberg2016optimal} that any other distributed control laws constructed by Laplacian matrices that do not correspond to complete graphs with equal edge weights are suboptimal.
\item When optimizing the edge weights, "competitions" are allowed between the connected agents. By doing so, the feasible region of the optimization problem is enlarged, and hence a smaller control energy cost might be obtained. We offer two sufficient conditions for when will "competitions" happen between agents. These two conditions help to determine whether the two agents will compete if we add a connection between them based on the old optimal solution.
\item When the graph topology is regular, namely, every node has the same number of neighbours, the controller can be calculated in a distributed manner.
\end{enumerate}

The rest of this paper is organized as follows. In Section \ref{sec:prelim}, some preliminaries and notations are presented. In Section \ref{sec:problem_formulation} we introduce the problem formulation. The design of the control gain is presented in Section \ref{sec:ctrl_g} and \ref{sec:opt_ctrl_gain}. In Section \ref{sec:opt_edge_weights}, the edge weights design is formulated as an SDP problem and we provide a different proof compared to \citep{thunberg2016optimal} from the angle of optimization and show that the lowest control energy cost is reached when the graph is complete and with equal edge weights. Two sufficient conditions on the existence of the negative optimal edge weight are presented. 
In Section \ref{sec:distri_opt}, a distributed way of computing the optimal edge weights is presented when the graph is regular. Finally, we conclude the paper and describe some future work in Section \ref{sec:conclusion}.

\section{Notations and Preliminaries}\label{sec:prelim}
In this section, we provide a brief introduction about some notations that we are going to use in the rest of the paper, as well as some basic knowledge about SDP. 

We denote $\bm{1}$ as an $N$ dimensional all-one column vector. $\bm{0}$ is denoted as an N dimensional all-zero matrix. The element on the $i$th row and $j$th column of any matrix $D$ is expressed as $[D]_{ij}$. $D_1\succeq D_2$ and $G_1\preceq G_2$ mean that $D_1-D_2$ and $G_2-G_1$ are positive semi-definite. $\otimes$ denotes the Kronecker product. $\|\cdot\|$ denotes 2-norm of matrice or vectors. $\|\cdot\|_F$ denotes the Frobenius norm of a matrix. We use $|\cdot|$ to denote the number of the elements of a set. And any notation with the superscript $*$ denotes the optimal solution to the corresponding optimization problem. $tr(\cdot)$ denotes the trace of a matrix. If $D,G\in\bm{S}^n_+$ are positive definite matrices, then $tr(DG)$ is the inner-product between $D$ and $G$. It also holds that $\nabla_Dtr(DG)=G$.

An edge-weighted undirected graph $\mathcal{G(V,E,W)}$ is composed of a node set $\mathcal{V}=\{1,2,\cdots,N\}$, an edge set $(i,j)\in\mathcal{E},\:i,j\in\mathcal{V}$ which describes the connection topology between the nodes and the edge weight set $w_{ij}\in\mathcal{W},\:i,j\in\mathcal{V}$ which includes all the weights of the corresponding edges. To abbreviate the notation, we label the edges with numbers. For example, an edge with label $l$ is denoted as $l\in\mathcal{E}$. On the other hand, seen from the nodes' perspective, the set of edges that is connected to node $i$ is denoted as $\mathcal{E}(i)$, which can be interpreted as communication channels of node $i$. Note that $\mathcal{E}=\bigcup_{i\in\mathcal{V}}\mathcal{E}(i)$. Similarly, the set of the edge weights belong to node $i$ is denoted as $\mathcal{W}(i)$. $\mathcal{N}(i)$ denotes the neighbour vertices set of node $i$.

Note that in this paper, we consider undirected graphs, hence the edge-weighted Laplacian matrix $L_w$ is symmetric and defined as
\begin{equation}\nonumber
	\begin{aligned}
		\left[L_w \right]_{ij}&=
		\begin{cases}
			\sum_lw_{il} &\mbox{if}\: i= j\:\mbox{and}\:(i,l)\in \mathcal{E}\\
			-w_{ij}&\mbox{if}\: i\neq j\:\mbox{and}\:(i,j)\in \mathcal{E}\\
			0 &\mbox{otherwise}
		\end{cases}\\
		\Leftrightarrow L_w&=\sum_{k\in \mathcal{E}}w_{k}E_{k},
	\end{aligned}
	\label{eq:Laplacian_matrix}
\end{equation}
where $k$ is the label of the edges, $w_k\in \mathcal{W},\:\forall k\in \mathcal{E}$ are the edge weights. If node $i$ and $j$ are connected via edge $k$, then
\begin{align}\nonumber
	\left[E_k\right]_{ii}=\left[E_k\right]_{jj}=1,\quad \left[E_k\right]_{ij}=\left[E_k\right]_{ji}=-1,
	\label{eq:E_k_def}
\end{align} 
and the other elements of $E_k$ are zero. For a connected graph, the eigenvalues of $L_w$ is denoted as $0=\lambda_1<\lambda_2\leq\cdots\leq\lambda_N$.

For any symmetric matrix $G\in\bm{S}^n$, $svec(G)$ is defined as
\begin{equation}\nonumber
	\begin{aligned}
	svec(G)=[G_{11},\sqrt{2}G_{21},\cdots,\sqrt{2}G_{n1},\\G_{22},\sqrt{2}G_{32},\cdots,\sqrt{2}G_{n1},\cdots,G_{nn}]^T.
	\end{aligned}
\end{equation}
It follows from the above definition that
\begin{equation}\nonumber
	tr(DG)=svec(D)^Tsvec(G),\quad \forall D,G\in \bm{S}^n.
	\label{eq:svec_tr}
\end{equation}
The \textit{Symmetric Kronecker Product} between two matrices $G$ and $D$ is defined by the following identity
\begin{equation}\nonumber
(R_1\otimes_s R_2)svec(G)=\frac{1}{2}svec(R_2GR_1^T+R_1GR_2^T),
\end{equation}
where $G\in\bm{S}^n$, but $R_1$ and $R_2$ is not necessarily symmetric. For more details about $svec(\cdot)$ and Symmetric Kronecker Product can be found in \cite{schacke2013kronecker} and \cite{alizadeh1998primal}.

Now we introduce some basics about SDP. An SDP problem in standard form can be expressed as
\begin{equation}
	\begin{aligned}
		& \underset{\xi_i}{\text{maximize}}
		& & \sum_{i=1}^m c_i\xi_i\\
		& \text{subject to}
		& &\sum_{i=1}^m\xi_iG_i\preceq R.
	\end{aligned}
	\label{eq:standard_sdp}
\end{equation}
Introducing the Lagrangian multiplier $\Phi\succeq 0$, the Lagrangian function of (\ref{eq:standard_sdp}) can be written as
\begin{equation}\nonumber
\mathcal{L}(\xi,\Phi)=-\sum_{i=1}^m c_i\xi_i+tr(\Phi(\sum_{i=1}^m \xi_iG_i-R))
\end{equation}
Taking the gradient of the Lagrangian function, we get its dual problem
\begin{equation}\nonumber
	\begin{aligned}
		& \underset{\Phi}{\text{minimize}}
		& & tr(R\Phi)\\
		& \text{subject to}
		& &tr(G_i\Phi)=c_i,\quad i=1,\cdots,m\\
		& & &\Phi\succeq 0.
	\end{aligned}
\end{equation}

\section{Problem Formulation}\label{sec:problem_formulation}
In this paper, we consider the following linear multi-agent system, where each agent has the following dynamics:
\begin{align}
	\dot{x}_i=Ax_i+Bu_i,\quad i=1,\cdots,N,
	\label{eq:system}
\end{align}
where $x_i\in \mathbb{R}^n, \forall i$ and they are connected through a edge-weighted graph $\mathcal{G(V,E,W)}$.

Suppose that $(A,B)$ is stabilizable, $A$ has no eigenvalue on the imaginary axis. We want each $x_i$ to reach consensus asymptotically, namely, $\|x_i(t)-x_j(t)\|\rightarrow 0, \:\forall i,j$ as $t\rightarrow\infty$, while minimizing the control energy cost. Plus, the control should be fully distributed and the control should only depend on relative information, i.e., the control must have the form
\begin{align}
	u_i=K\sum_{j\in \mathcal{N}(i)}w_{ij} (x_i-x_j).
	\label{eq:control}
\end{align}
Note here that, compared to the classical definition, $w_{ij}$ does not have to be positive, meaning that we allow some of the agents to compete against each other. The reason for us to allow the edge weights to be negative will be described later.

Now we formulate the problem as an infinite horizon Linear Quadratic (LQ) optimal control problem:
\begin{equation}
	\begin{aligned}
		& \underset{K,\{w_{ij}\}}{\text{minimize}}
		& & J(K,\{w_{ij}\})= \int_0^\infty U^TU\:dt\\
		& \text{subject to}
		& &\dot{X}=(I_N\otimes A)X+(I_N\otimes B)U,\\
		& & &U=(L_w\otimes K)X,\\
		& & &\|x_i(t)-x_j(t)\|\rightarrow 0,\:t\rightarrow\infty, \:\forall i,j,
	\end{aligned}
	\tag{LQ1}\label{eq:LQR_pro}
\end{equation}
where $X=[x_1^T,\cdots,x_N^T]^T$, $U=[u_1^T,\cdots, u_N^T]^T$, and $L_w$ is the weighted Laplacian matrix of the edge-weighted graph $\mathcal{G(V,E,W)}$.

The problem (\ref{eq:LQR_pro}) is well-posed though $\lim_{t\rightarrow\infty}x_i(t)=\infty$ when $A$ is not stable, namely, $J(K^*,\{w_{ij}^*\})<\infty$. To see this, we simplify the system in (\ref{eq:LQR_pro}) by applying a coordinate change  $\tilde{X}=(T^T\otimes I_n)X$, where $T^TT=I$ ,$T^TL_wT=\Lambda_L$ and $\Lambda_L=diag(0,\lambda_2,\cdots,\lambda_N)$, \citep{fax2004information}. Then the closed-loop system becomes:
\begin{align}\nonumber
	\dot{\tilde{X}}=(I_N\otimes A)\tilde{X}+(\Lambda_L\otimes BK)\tilde{X}.
\end{align}
Note that if the systems (\ref{eq:system}) reach consensus, then $\lim_{t\rightarrow\infty}X(t)=\bm{1}\otimes\bar{X}(t)$, where $\bar{X}:\mathbb{R}\mapsto \mathbb{R}^n$. $T$ diagonalizes the edge-weighted Laplacian matrix $L_w$, hence it should have that $T=[\frac{1}{N}\bm{1},t_2,\cdots,t_N]$, where $t_i$ is the $i$th column of matrix $T$ and $t_i\perp\bm{1}$. It follows that
\begin{equation}\nonumber
\begin{aligned}
\lim_{t\rightarrow\infty}\tilde{X}(t)&=\lim_{t\rightarrow\infty}(T^T\otimes I_n)X(t)\\
&=(
\begin{pmatrix}
\frac{1}{N}\bm{1}^T\\t_2^T\\\vdots\\t_N^T
\end{pmatrix}\otimes I_n)(\bm{1}\otimes\bar{X}(t))
=\begin{pmatrix}
1\\0\\\vdots\\0
\end{pmatrix}\otimes \bar{X}(t)
\end{aligned}
\end{equation}
and therefore $\lim_{t\rightarrow\infty}\tilde{x}_i(t)=0,i\geq 2$ if consensus is reached in the original system (\ref{eq:system}). On the other hand, if $\lim_{t\rightarrow\infty}\tilde{x}_i(t)=0,i\geq 2$, then $\lim_{t\rightarrow\infty}X(t)=\lim_{t\rightarrow\infty}(T\otimes I_n)\tilde{X}(t)=\lim_{t\rightarrow\infty}([\frac{1}{N}\bm{1},t_2,\cdots,t_N]\otimes I_n)([1,0 \cdots,0]^T\otimes\tilde{x}_1(t))=\frac{1}{N}\bm{1}\otimes\tilde{x}_1(t)$. Therefore, letting (\ref{eq:system}) reach asymptotic consensus is equivalent to regulate $\tilde{x}_i,\:i=2,\cdots,N$. Furthermore, the cost function $J$ in (\ref{eq:LQR_pro}) reads:
\begin{equation}
	\begin{aligned}
		&J(K,\{w_{ij}\})=\int_0^\infty U^TU\:dt=\int_0^\infty X^T(L_w^2\otimes K^TK)X\:dt\\
		&=\int_0^\infty \tilde{X}^T(T^T\otimes I_n)(L_w^2\otimes K^TK)(T\otimes I_n)\tilde{X}\:dt\\
		&=\int_0^\infty \tilde{X}^T(T^TL_w\underbrace{TT^T}_IL_wT\otimes K^TK)\tilde{X}\:dt\\
		&=\int_0^\infty \tilde{X}^T(\Lambda_L^2\otimes K^TK)\tilde{X}\:dt.\\
	\end{aligned}
	\label{eq:tilde_cost_function}
\end{equation}
Note that $\Lambda_L=diag(0,\lambda_2,\cdots,\lambda_N)$, hence the $\tilde{x}_1$ does not appear in the cost function $J$ and we can remove everything about $\tilde{x}_1$ and simplify the problem (\ref{eq:LQR_pro}) as
\begin{equation}
	\begin{aligned}
		& \underset{K,\{w_{ij}\}}{\text{minimize}}
		& & J(K,\{w_{ij}\})= \int_0^\infty \hat{U}^T\hat{U}\:dt\\
		& \text{subject to}
		& &\dot{\hat{X}}=(I_{N-1}\otimes A)\hat{X}+(I_{N-1}\otimes B)\hat{U},\\
		& & &\hat{U}=(\Lambda_L^\prime\otimes K)\hat{X},
	\end{aligned}
	\tag{LQ2}\label{eq:LQR_pro1}
\end{equation}
where $\Lambda_L^\prime=diag(\lambda_2,\cdots,\lambda_N)$ and $\hat{X}=[\tilde{x}_2^T,\cdots,\tilde{x}_N^T]^T$. $\Lambda_L^\prime$ is determined by the edge weights $w_{ij}$.

\section{Main Results}
The product between $\Lambda_L^\prime$ and $K$ makes (\ref{eq:LQR_pro1}) a non-convex problem. Moreover, even when $\Lambda_L^\prime$ is given, similar to \citep{borrelli2008distributed}, computing such control gain $K$ while minimizing the energy cost is still hard. Hence we focus on the control gain $K$ that guarantees asymptotic consensus first and try to find an optimal control among this family of controllers.

\subsection{The Control That Guarantees Consensus}\label{sec:ctrl_g}
Consider the following control gain
\begin{align}
	K=-\frac{1}{\lambda_2}B^TP,
	\label{eq:ctrl_design}
\end{align}
where $P$ is the unique positive semi-definite stabilizing solution (see Sec 10.3, \citep{trentelman2012control}) to the following Algebraic Riccati Equation (ARE):
\begin{align}
	A^TP+PA-PBB^TP=-Q,
	\label{eq:ARE}
\end{align}
where $Q\succeq 0$.
Note that the control design (\ref{eq:ctrl_design}) is related to the results in \citep{zhang2011optimal}, but they require $Q\succ 0$ while in our case we only require $Q\succeq 0$. Hence the set of $Q$ that we can choose from is closed.

\begin{prop}
	The control (\ref{eq:ctrl_design}) will let the multi-agent system (\ref{eq:system}) reach consensus asymptotically if the weighted Laplacian matrix $L_w$ is positive semi-definite and has only one zero eigenvalue.
\end{prop}
\begin{pf}
Since $P$ is the unique stabilizing solution to the ARE (\ref{eq:ARE}), then $A_p=A-BB^TP$ must be Hurwitz. If we plug in the control gain (\ref{eq:ctrl_design}), the closed-loop system becomes $\dot{\tilde{x}}_i=(A-\frac{\lambda_i}{\lambda_2}BB^TP)\tilde{x}_i=(A-\sigma_iBB^TP)\tilde{x}_i=A_{cl}^i\tilde{x}_i$, where $\sigma_i=\frac{\lambda_i}{\lambda_2},\: i=2,\cdots,N$. $L_w$ is positive semi-definite and has only one zero eigenvalue and hence $\sigma_i\geq 1$.

By rewriting (\ref{eq:ARE}), we get
\begin{align}
&A^TP+PA-PBB^TP=-Q\nonumber\\
\Leftrightarrow &A_{cl}^{iT}P+PA_{cl}^i=-Q-(2\sigma_i-1)PBB^TP.\label{eq:ARE_lyap}
\end{align}
Suppose $(\lambda_{cl},v_{cl})$ is an eigen-pair of $A_{cl}^i$, then multiply $v_{cl}^H$ on the left and $v_{cl}$ on the right on both hand sides of (\ref{eq:ARE_lyap}), it follows that
\begin{equation}
2Re(\lambda_{cl})\cdot v_{cl}^HPv_{cl}=-v_{cl}^HQv_{cl}-(2\sigma_i-1)v_{cl}^HPBB^TPv_{cl},
\label{eq:re_neg}
\end{equation}
where $v_{cl}^H$ denotes the conjugate transpose of $v_{cl}$.
The eigen-pair $(\lambda_{cl},v_{cl})$ must satisfy one of the following three cases:
\begin{enumerate}
\item If $v_{cl}^HPv_{cl}=0$, then $Pv_{cl}=0$. Hence $A_{cl}^iv_{cl}=(A-\sigma_iBB^TP)v_{cl}=Av_{cl}=\lambda_{cl}v_{cl}$. On the other hand, $A_pv_{cl}=(A-BB^TP)v_{cl}=Av_{cl}=\lambda_{cl}v_{cl}$, thus $(\lambda_{cl},v_{cl})$ is also an eigen-pair of $A_p$. We know that $A_p$ is Hurwitz, and hence $Re(\lambda_{cl})<0$.
\item If $v_{cl}^HPv_{cl}>0$ and right hand side of (\ref{eq:re_neg}) is strictly negative, then it is easy to conclude that $Re(\lambda_{cl})<0$.
\item If $v_{cl}^HPv_{cl}>0$ and right hand side of (\ref{eq:re_neg}) is zero, then it follows that $Re(\lambda_{cl})=0$. It also follows that $v_{cl}^HQv_{cl}=0$ and $v_{cl}^HPBB^TPv_{cl}=0$ since both $Q$ and $PBB^TP$ is positive semi-definite. Hence $B^TPv_{cl}=0$. Then we can also conclude that $A_{cl}^iv_{cl}=Av_{cl}=A_pv_{cl}=\lambda_{cl}v_{cl}$ and hence $(\lambda_{cl},v_{cl})$ is also an eigen-pair of $A_p$. We know that $A_p$ is Hurwitz, thus $Re(\lambda_{cl})<0$. This contradicts $Re(\lambda_{cl})=0$ derived from (\ref{eq:re_neg}). Thus this case will never happen.
\end{enumerate}

Since any eigenvalue $\lambda_{cl}$ have $Re(\lambda_{cl})<0$, then $A_{cl}^i$ is Hurwitz. Therefore control (\ref{eq:ctrl_design}) will let the multi-agent system (\ref{eq:system}) reach consensus asymptotically.\qed
\end{pf}
Note that the control gain (\ref{eq:ctrl_design}) depends on the choice of $Q$, hence we would like to choose a good $Q$ so that the control gain $K$ should optimize the energy cost. 

\subsection{Interval Bound on the Energy Cost}\label{sec:opt_ctrl_gain}
Plugging the dynamics of the closed-loop system and the control gain (\ref{eq:ctrl_design}) into the cost function $J(K,\{w_{ij}\})$, we get
\begin{equation}
\begin{aligned}
&J(K,\{w_{ij}\})=\int_0^\infty \tilde{X}^T(\Lambda_L^2\otimes K^TK)\tilde{X}\:dt\\
&=\int_0^\infty \hat{X}^T(\Lambda_L^{\prime 2}\otimes \frac{1}{\lambda_2^2}PBB^TP)\hat{X}\:dt=\hat{X}_0^TH\hat{X}_0,
\end{aligned}
\end{equation}
where $H=\int_0^\infty e^{(I_{N-1}\otimes A-\Sigma_L^\prime\otimes BB^TP)^Tt}(\Sigma_L^{\prime 2}\otimes PBB^TP)$ $e^{(I_{N-1}\otimes A-\Sigma_L^\prime\otimes BB^TP)t}\:dt$, $\Sigma_L^\prime=diag(\sigma_2,\cdots,\sigma_N)$ and $\sigma_i=\lambda_i/\lambda_2$. Note that $H$ is block diagonal, i.e., $H=diag(H_2,\cdots,H_N)$, where
\begin{equation}
H_i=\int_0^\infty e^{A_{cl}^{iT}t}\sigma_i^2PBB^TPe^{A_{cl}^it}\:dt.
\label{eq:H_i}
\end{equation}

\begin{thm}\label{thm:bound on H_i}
$P_0\preceq H_i\preceq \frac{\sigma_i^2}{2\sigma_i-1}P$, where $P_0$ is the solution to the ARE (\ref{eq:ARE}) when $Q=0$ and $P$ is the solution to the ARE (\ref{eq:ARE}) for any $Q\succeq 0$.
\end{thm}
\begin{pf}
We first show the upper bound. 
Note that $H_i$ has the form of (\ref{eq:H_i}), hence it is actually the analytic solution to the following Lyapunov equation
\begin{equation}
A_{cl}^{iT}H_i+H_iA_{cl}^i=-\sigma_i^2PBB^TP.
\label{eq:H_i lyap}
\end{equation}
Multiply $\frac{\sigma_i^2}{2\sigma_i-1}$ on both hand sides of (\ref{eq:ARE_lyap}) and subtract (\ref{eq:H_i lyap}), we have
\begin{equation}
\begin{aligned}
A_{cl}^{iT}(\frac{\sigma_i^2}{2\sigma_i-1}P-H_i)+(\frac{\sigma_i^2}{2\sigma_i-1}P-H_i)A_{cl}^i\\
=-\frac{\sigma_i^2}{2\sigma_i-1}Q,
\end{aligned}\nonumber
\end{equation}
which is a Lyapunov equation. Consequently, since $Q\succeq 0$ and $A_{cl}^i$ is stable, it must hold that 
\begin{align}
\frac{\sigma_i^2}{2\sigma_i-1}P-H_i=\int_0^\infty e^{A_{cl}^{iT}t}\frac{\sigma_i^2}{2\sigma_i-1}Qe^{A_{cl}^it}\:dt\succeq 0.\label{eq:analytic sol}
\end{align}
Hence $\frac{\sigma_i^2}{2\sigma_i-1}P\succeq H_i$.

Now we show the lower bound. When $Q=0$, the solution $P_0$ must satisfy
\begin{align}
&A^TP_0+P_0A-P_0BB^TP_0=0\nonumber\\
\Leftrightarrow&A_{cl}^{iT}P_0+P_0A_{cl}^i=-\sigma_iPBB^TP_0-\sigma_iP_0BB^TP\nonumber\\&+P_0BB^TP_0.\label{eq:ARE_P0 rewrite}
\end{align}
Subtract (\ref{eq:H_i lyap}) by (\ref{eq:ARE_P0 rewrite}), we get
\begin{equation}\nonumber
\begin{aligned}
A_{cl}^{iT}(H_i-P_0)+(H_i-P_0)A_{cl}^i\\=-(\sigma_iP-P_0)BB^T(\sigma_iP-P_0).
\end{aligned}
\end{equation}
Also, since $(\sigma_iP-P_0)BB^T(\sigma_iP-P_0)\succeq 0$ and $A_{cl}^i$ is Hurwitz, it must hold that 
\begin{align}
H_i-P_0&=\int_0^\infty e^{A_{cl}^{iT}t}(\sigma_iP-P_0)BB^T\nonumber(\sigma_iP-P_0)e^{A_{cl}^it}\:dt\nonumber\\&\succeq 0.\nonumber
\end{align}
Hence $H_i\succeq P_0$.
\qed
\end{pf}
Note that the interval bound will degenerate to a "point" if $Q=0$ is chosen and when the graph is complete, namely, in this case $P=P_0$ and $\sigma_i=1,\:i=2,\cdots,N$. Since by definition $\sigma_i\geq 1$, then $\frac{\sigma_i^2}{2\sigma_i-1}$ is monotonously increasing with respect to $\sigma_i$, hence $P_0\preceq H_i\preceq \frac{\sigma_i^2}{2\sigma_i-1}P\preceq \frac{\sigma_N^2}{2\sigma_N-1}P$.

Because the optimal control energy cost has a lower bound, the distance between the optimal control energy cost and the upper bound is also bounded and hence it is a rational choice to optimize the upper bound,i.e., $\frac{\sigma_N^2}{2\sigma_N-1}P$,  since it is hard to optimize the real energy cost. Note that the upper bound consists of two variables: $\sigma_N$ and $P$ and hence is optimized from two perspectives. From the perspective of $P$, based on the monotonicity of the solution to the ARE (\cite{willems1971least}) (namely, if $Q_1\succeq Q_2$, then the corresponding solution holds $P_1\succeq P_2$), the upper bound is minimized by choosing $Q=0$. More precisely, when $Q=0$, the right hand side of (\ref{eq:analytic sol}) becomes zero and hence $H_i=\frac{\sigma_i^2}{2\sigma_i-1}P_0$.

Now that we have optimized the upper bound from the perspective of $P$, what remains to optimize is $\sigma_N$, which is determined by the edge weights. This will be done in the next section.

\subsection{Optimizing the Edge Weights}\label{sec:opt_edge_weights}
$\frac{\sigma_N^2}{2\sigma_N-1}$ is monotonously increasing since $\sigma_N\geq1$, hence minimizing $\frac{\sigma_N^2}{2\sigma_N-1}$ is equivalent to minimizing $\sigma_N$.
\begin{rem}
	$\sigma_N=\frac{\lambda_N}{\lambda_2}$ is so-called ``synchronizability'' in the network synchronization field on physics, and its optimization has been widely studied in the past. \cite{jalili2013enhancing} gives a thorough survey on this topic. Though there has been a lot of approaches to enhance synchronizability in the field of physics \citep{jalili2013enhancing}, heuristic approaches are used to adjust the edge weights and the effectiveness of the approaches are shown empirically by numerical examples. We would like to know what is the exact optimal edge weight realization to minimize the synchronizability. Since the problem involves the eigenvalues of the matrix $L_w$, it is natural to formulate the problem as an SDP problem.\qed
\end{rem}

Recall that in the problem formulation, $L_w$ is symmetric and negative edge weights $w_{ij}$ are allowed. Furthermore, it is assumed that there is a finite communication resource in the network, namely, $\sum_{k\in\mathcal{E}}w_k^2=const$. Without loosing generality, we assume the constant to be 1. Of course, the optimized edge-weighted graph $\mathcal{G}$ should be connected, namely, the second smallest eigenvalue of the Laplacian matrix $L_w^*$ should be strictly positive. To ensure the existence of feasible edge weight realization, we make the assumption that the unweighted graph represented by the unweighted Laplacian matrix $L_u=\sum_{k\in\mathcal{E}}E_k$ is connected. Thus, the problem can be formulated as the following optimization problem:
\begin{equation}
	\begin{aligned}
		& \underset{\lambda_2,\lambda_N,\mu,w_k}{\text{minimize}}
		& & \frac{\lambda_N}{\lambda_2} \\
		& \text{subject to}
		& & \lambda_2I-\mu\bm{1}\bm{1}^T\preceq \sum_{k=1}^{|\mathcal{E}|} w_kE_k\preceq \lambda_NI,\\
		& & & \sum_{k=1}^{|\mathcal{E}|}w_k^2=1,\\
		& & & \lambda_2>0.
	\end{aligned}
	\tag{P1}\label{eq:convex_opt_edge_add1}
\end{equation}
In (\ref{eq:convex_opt_edge_add1}), the variable $\mu$ is used to shift the zero eigenvalue of the Laplacian matrix $L_w$ with its eigenvector $\bm{1}$. When the optimal value is reached, $\lambda_2^*$ would be the smallest eigenvalue of $\sum_{k=1}^{|\mathcal{E}|} w_kE_k+\mu\bm{1}\bm{1}^T$. By using the property that for any semi-positive definite matrix $G$, we have $\beta_1I\preceq G\preceq\beta_nI$, where $\beta_1,\beta_n$ is the smallest and biggest eigenvalue of $G$ respectively, we get the above constraints. 
\begin{rem}
The reason we put a free variable $\mu$ in the constraint $\lambda_2I-\mu\bm{1}\bm{1}^T\preceq \sum_{k=1}^{|\mathcal{E}|} w_kE_k$ instead of writing the constraint as $\lambda_2(I-\frac{1}{N}\bm{1}\bm{1}^T)\preceq \sum_{k=1}^{|\mathcal{E}|} w_kE_k$ is that the latter formulation does not have any interior point in the feasible domain, i.e., every feasible point lies on the boundary. Later we will consider the dual problem and we want to use the Slater's condition to show the strong duality holds, which requires the primal problem to have at least one strictly feasible solution. Note that when optimal value is reached, $\mu^*$ does not have to be $\lambda_2^*/N$. In fact, given an optimal edge weight realization and $\lambda_2^*$, any $\mu\geq \lambda_2^*/N$ is an optimal solution. 
\end{rem}

\begin{rem}
Note that in (\ref{eq:convex_opt_edge_add1}), we allow the agents to "compete" by not introducing the constraints $w_k\geq 0,\:k=1,\cdots,|\mathcal{E}|$ into the problem. If we added the constraints of non-negative edge weights, the feasible domain would shrink and hence potentially gives a larger optimal value. Thus, allowing the agents to compete will give an at least as good, sometimes better system performance, which is quite contrast to intuition.
\qed
\end{rem}
Though the problem (\ref{eq:convex_opt_edge_add1}) is not convex, it can be rewritten into a convex problem. By applying variable change $t=\lambda_N/\lambda_2$, $\gamma=1/\lambda_2$, $\mu/\lambda_2=y_0$, $w_k/\lambda_2=y_k$ and $1/\lambda_2=\gamma$, we change the problem into a semidefinite programming (SDP) problem:
\begin{equation}
	\begin{aligned}
		& \underset{\gamma,t,y_k,y_0}{\text{minimize}}
		& & t\\
		& \text{subject to}
		& & I-y_0\bm{1}\bm{1}^T\preceq\sum_{k=1}^{|\mathcal{E}|}y_kE_k\preceq tI,\\
		& & &\gamma>0.
	\end{aligned}
	\tag{P2}\label{eq:cvx_pro}
\end{equation}


Note that the last constraint $\gamma>0$ in \eqref{eq:cvx_pro} is redundant since it does not affect the other constraints and the objective function. We can just remove the constraint and get a classic SDP problem. It is also worth noticing that with the aforementioned change of variables, the control \eqref{eq:control} converges to the following form:
\begin{align*}
	u_i
	=-B^TP_0\sum_{j\in\mathcal{N}(i)}y_{ij}(x_i-x_j).
\end{align*}
This means that $J(-\frac{1}{\lambda_2}B^TP_0,\{w_{ij}\})=J(-B^TP_0,\{y_{ij}\})$. (It depends on how we denote $K$.)

We know that the optimal value of the primal problem can always be attained. Hence if we are able to show that the Slater's condition holds for (\ref{eq:cvx_pro}), i.e., there exists a solution $(t,y_k,y_0)$ which is strictly feasible, then the strong duality holds. Note that for the constraint $\sum_{k=1}^{|\mathcal{E}|}y_kE_k\preceq tI$, we can always choose $t$ to be big enough such that the strict inequality holds. For the first matrix constraint, $w_k=\frac{1}{\lambda_2^*-\varepsilon}$, $\mu=1+\varepsilon$ is a strictly feasible solution given that $\lambda_2^* >\varepsilon>0$ \citep{goring2008embedded}. Hence the Slater's condition holds and strong duality holds between (\ref{eq:cvx_pro}) and (\ref{eq:dual_pro_simp}).

By Lagrangian duality, one can get the dual problem of (\ref{eq:cvx_pro}).
\begin{equation}
	\begin{aligned}
		& \underset{\Phi_1,\Phi_2}{\text{maximize}}
		& & tr(\Phi_1)\\
		& \text{subject to}
		& & tr(\Phi_2)=1,\\
		& & & tr(\bm{1}\bm{1}^T\Phi_1)=0,\\
		& & & tr(E_k\Phi_2)-tr(E_k\Phi_1)=0,\: k=1,2,\cdots,|\mathcal{E}|,\\
		& & & \Phi_1,\Phi_2\succeq 0.\\
	\end{aligned}
	\tag{D2}\label{eq:dual_pro_simp}
\end{equation}
The Karush-Kuhn-Tucker (KKT) conditions for (\ref{eq:cvx_pro}) reads:
\begin{equation}
	\begin{aligned}
		&tr(\Phi_2^*)=1,\quad tr(\bm{1}\bm{1}^T\Phi_1^*)=0,\\
		&tr(E_k\Phi_2^*)-tr(E_k\Phi_1^*)=0,\:k=1,\cdots,|\mathcal{E}|,\\
		&tr(\Phi_1^*)=\sum_{k=1}^{|\mathcal{E}|}y_k^*tr(E_k\Phi_1^*),\quad t^*=\sum_{k=1}^{|\mathcal{E}|}y_k^*tr(E_k\Phi_2^*),
	\end{aligned}
	\label{eq:KKT_cond}
\end{equation}
with the primal feasiblity constraints ommitted here for the sake of brevity.

\cite{thunberg2016optimal} considered a "topology free" problem, namely, their problem formulation is similar to (\ref{eq:LQR_pro}), but without the constraint $U=(L_w\otimes K) X$. Their control is distributed and turns out to have the form of $U=(L_w\otimes K) X$ as well, whose edge-weighted Laplacian matrix in the controller corresponds to a complete graph with equal weights $\frac{1}{N}$ on each edge. Thus their controller is optimal to any other control with a different control gain, graph topologies which is not complete and edge weights not equal to $\frac{1}{N}$.
However, \cite{thunberg2016optimal} have not shown the uniqueness of the optimal solution (which can actually be shown using the uniqueness of the solution for the ARE). Here we provide another proof from the view of an optimization problem. Plus, we show the uniqueness of the optimal solution.

\begin{thm}[c.f. \cite{thunberg2016optimal}]
Among all controls that belong to the family \eqref{eq:ctrl_design}, if given the freedom of choosing positive semidefinite matrix $Q$ and the graph topology among $N$ nodes as long as its weighted Laplacian matrix $L_w\succeq 0$, the optimal energy cost is reached by the control generated by $Q=0$ and $L_w$ corresponding to a complete graph with equal edge weights.
\end{thm}
\begin{pf}
	In this case, the matrices $E_k$ in all should represent the topology of a fully connected graph so that (\ref{eq:cvx_pro}) can have the biggest feasible domain.
	Choose $y_k^*=\frac{1}{N},\:k=1,\cdots,|\mathcal{E}|$, $t^*=1$ and
	\begin{equation}
		\Phi_1^*=\Phi_2^*=
		\begin{bmatrix}
			\frac{1}{N} &-\frac{1}{N(N-1)} &\cdots &-\frac{1}{N(N-1)}\\
			-\frac{1}{N(N-1)} &\frac{1}{N} &\ddots &\vdots\\
			\vdots &\cdots &\ddots &-\frac{1}{N(N-1)}\\
			-\frac{1}{N(N-1)} &\cdots &-\frac{1}{N(N-1)} &\frac{1}{N}
		\end{bmatrix},
		\label{eq:Phi1Phi2}
	\end{equation}
	which is feasible to the dual problem (\ref{eq:dual_pro_simp}). All KKT conditions are then trivially fulfuilled, except for the complementarity slackness.
	Plugging in $t^*=1$ and $y_k^*=\frac{1}{N}$, the complementarity slack conditions read as
	\begin{align}
		&tr\lbrack(I-y_0^*\bm{1}\bm{1}^T-\sum_{k=1}^{\frac{n(n-1)}{2}}y_k^*E_k)\Phi_1^*\rbrack\nonumber\\
		&=1-\frac{1}{N}\sum_{k=1}^{\frac{n(n-1)}{2}}tr(E_k\Phi_1^*)
		=1-\frac{1}{N} tr(L_u^f\Phi_1^*)=0,
		\label{eq:compl_slack1}
	\end{align}
	\begin{align}
		tr\lbrack(\sum_{k=1}^{\frac{n(n-1)}{2}}y_k^*E_k-t^*I)\Phi_2^*\rbrack&=
		\frac{1}{N}\sum_{k=1}^{\frac{n(n-1)}{2}}tr(E_k\Phi_2^*)-1\nonumber\\&=\frac{1}{N} tr(L_u^f\Phi_2^*)-1=0,
		\label{eq:compl_slack2}
	\end{align}
	where $L_u^f$ denotes the Laplacian matrix of a fully-connected unweighted graph. 
Since $\Phi_1^*=\Phi_2^*$, (\ref{eq:compl_slack1}) is the same as (\ref{eq:compl_slack2}). And we know that $tr(\Phi_1^*)=1$ and $tr(\bm{1}\bm{1}^T\Phi_1^*)=0$, then $\frac{1}{N} tr\lbrack (L_u^f+\bm{1}\bm{1}^T)\Phi_1^*\rbrack=\frac{1}{N}Ntr(\Phi_1^*)=1$. Therefore, the KKT conditions are satisfied and $y_k^*=\frac{1}{N}$ is an optimal solution to the problem (\ref{eq:cvx_pro}).

	Now we show the uniqueness of the solution. Suppose apart from $y_k^*=\frac{1}{N},\forall k$, there is another optimal solution $y_k^{*\prime}$ also reaches the optimal value $t^*=1$. Both $y_k^*$ and $y_k^{*\prime}$ would have to satisfy the constraints, namely:
	\begin{align}
		I-y_0^*\bm{1}\bm{1}^T\preceq \sum_{k=1}^{\frac{n(n-1)}{2}}y_k^{*}E_k\preceq I\label{eq:opt_constraint1}\\
		I-y_0^{*\prime}\bm{1}\bm{1}^T\preceq \sum_{k=1}^{\frac{n(n-1)}{2}}y_k^{*\prime}E_k\preceq I\label{eq:opt_constraint2}
	\end{align}
	By subtracting (\ref{eq:opt_constraint1}) and (\ref{eq:opt_constraint2}), we have
	\begin{equation}
	\begin{aligned}
		\sum_{k=1}^{\frac{N(N-1)}{2}}(y_k^*-y_k^{*\prime})E_k\preceq y_0^{*\prime}\bm{1}\bm{1}^T,\nonumber\\
		\sum_{k=1}^{\frac{N(N-1)}{2}}(y_k^{*\prime}-y_k^*)E_k\preceq y_0^*\bm{1}\bm{1}^T.
	\end{aligned}
	\end{equation}
	
	This means for any $v\perp\bm{1}$, it must hold
	\begin{equation}
	\begin{aligned}
		\sum_{k=1}^{\frac{N(N-1)}{2}}(y_k^*-y_k^{*\prime})v^TE_kv\leq 0,\nonumber\\
		\sum_{k=1}^{\frac{N(N-1)}{2}}(y_k^{*\prime}-y_k^*)v^TE_kv\leq 0.
	\end{aligned}
	\end{equation}
	Since $E_k\succeq 0,\:\forall k$, then it must hold $y_k^*=y_k^{*\prime}$.

Note that since the optimal value of the problem (\ref{eq:convex_opt_edge_add1}) is 1, this means that $\lambda_2=\lambda_3=\cdots=\lambda_N$. Recall that the energy cost function is $J(-\frac{1}{\lambda_2}B^TP,\{\frac{1}{N}\})=\hat{X}_0^TH\hat{X}_0$ and hence the bounds in Theorem \ref{thm:bound on H_i} degenerate to $P_0\preceq H_i\preceq P,\:i=2,\cdots,N$. If we choose $Q=0$, then every block $H_i=P_0$ and the optimal energy cost control is obtained. Thus we prove the theorem.\qed
\end{pf}

The above theorem also describes the fact that our control design will coincides with the controller in \cite{thunberg2016optimal} when the graph is complete. We illustrate the above theorem by the following example. 

\begin{exmp}
Consider a complete graph with 3 nodes.
	\begin{figure}[!htpb]
		\centering
		\includegraphics[width=0.3\textwidth]{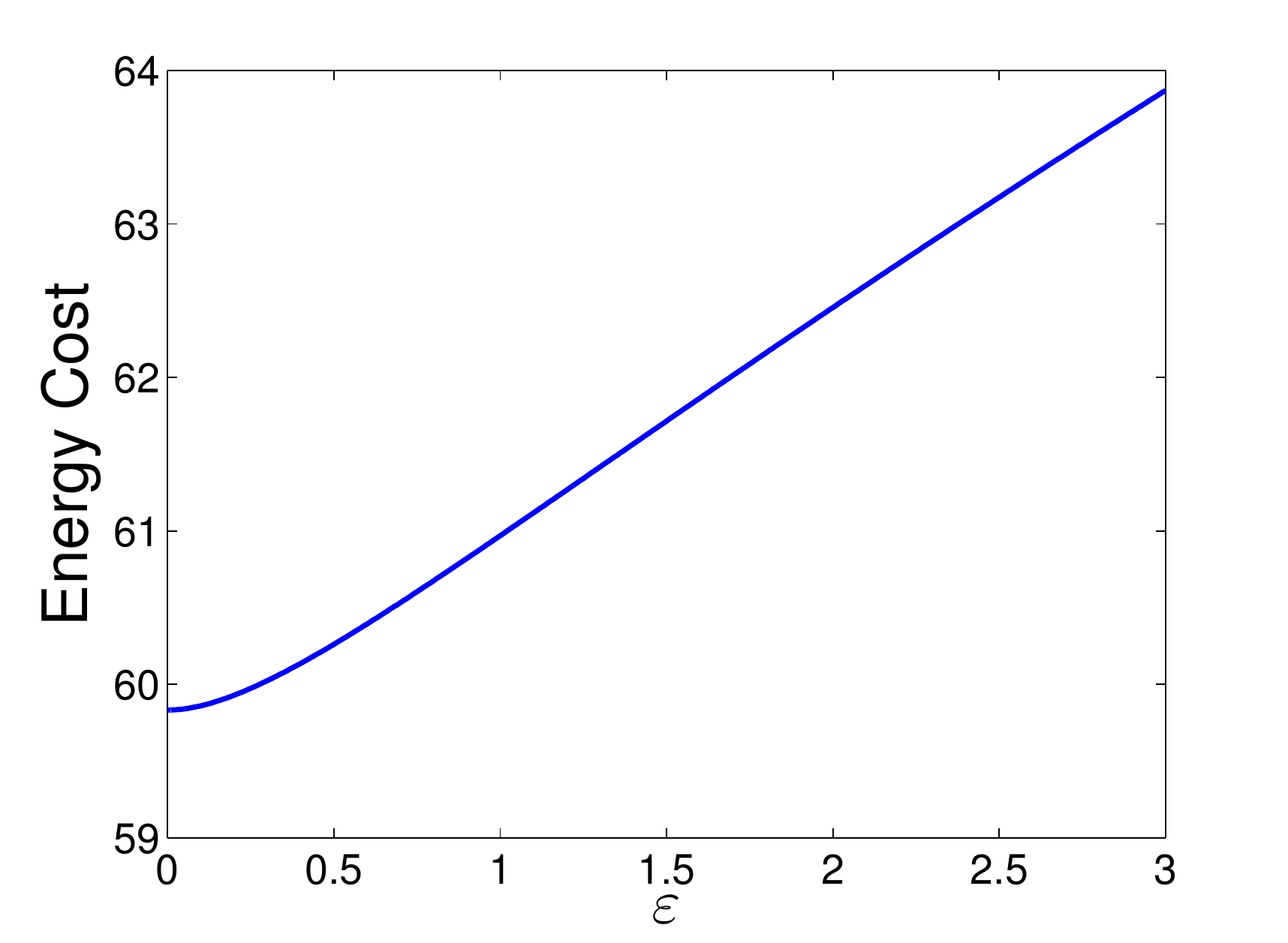}
		\caption{Energy Cost VS. $\varepsilon$.}
		\label{fig:energy_cost_epsilon}
	\end{figure}
\end{exmp}
We choose the matrices in \eqref{eq:system} and ARE as
\begin{equation}
A=
\begin{bmatrix}
0 &1\\-1 &2
\end{bmatrix},\quad
B=
\begin{bmatrix}
1\\2
\end{bmatrix},\quad
Q=\varepsilon I,
\label{eq:A B Q}
\end{equation}
and $\varepsilon$ varies from zero to three. The initial value $\hat{X}_0$ is chosen as $[-1.3077, -0.4336,0.3426,3.5784]^T$ and the same initial value is used when $\varepsilon$ varies. We can see from Fig. \ref{fig:energy_cost_epsilon} that when $\varepsilon=0$, the energy cost is minimum; and when $\varepsilon$ increases, the energy cost also increases.

An interesting question to ask is: given an existing optimal edge weight realization of a graph, if a new edge is added between two nodes, when can an optimal weight of the newly added edge be strictly negative? Recall that the control (\ref{eq:control}), when $w_{ij}<0$, this means that the two connected agents $i$ and $j$ are ``competing with'' or ``pushing away from'' each other but the entire system performance is still optimal. We would like to understand when will this happen. Here we provide two sufficient conditions for the existence of such optimal solutions.
\begin{thm}
	Suppose $\Phi_1^*$ and $\Phi_2^*$ are the optimal solutions to the dual problem (\ref{eq:dual_pro_simp}) for the graph $\mathcal{G(V,E,W)}$. Denote $\mathcal{\hat{G}(V,\hat{E},\hat{W})}$, where $\mathcal{\hat{E}}=\mathcal{E}\bigcup\{e_{|\mathcal{E}|+1}\}$, $e_{|\mathcal{E}|+1}\not\in \mathcal{E}$, $\mathcal{\hat{W}}=\mathcal{W}\bigcup\{w_{|\mathcal{E}|+1}\}$, $w_{|\mathcal{E}|+1}\not\in \mathcal{W}$. If $tr(E_{|\mathcal{E}|+1}\Phi_2^*)-tr(E_{|\mathcal{E}|+1}\Phi_1^*)=0$, then there must exist at least one optimal edge-weight realization of $\mathcal{\hat{G}(V,\hat{E},\hat{W})}$ where $w_{|\mathcal{E}|+1}^*<0$ and $t^*=\hat{t}^*$.
\end{thm}
\begin{pf}
Note that $\mathcal{\hat{G}(V,\hat{E},\hat{W})}$ is adding one constraint in the dual. If $tr(E_{|\mathcal{E}|+1}\Phi_2^*)-tr(E_{|\mathcal{E}|+1}\Phi_1^*)=0$, then the optimal solution $\Phi_1^*$ and $\Phi_2^*$ of the dual problem corresponding to $\mathcal{G(V,E,W)}$ should also be an optimal solution to the dual problem corresponding to $\mathcal{\hat{G}(V,\hat{E},\hat{W})}$.

	Together with the KKT conditions for the optimal edge-weight realization of $\mathcal{G(V,E,W)}$, we have the following:
	\begin{equation}
		\begin{aligned}
			&tr(\Phi_2^*)=1,\quad tr(\bm{1}\bm{1}^T\Phi_1^*)=0,\\
			&tr(E_k\Phi_2^*)-tr(E_k\Phi_1^*)=0,\:k=1,\cdots,|\mathcal{E}|,\\
			&tr(E_{|\mathcal{E}|+1}\Phi_2^*)-tr(E_{|\mathcal{E}|+1}\Phi_1^*)=0,\\
			&tr(\Phi_1^*)=\sum_{k=1}^{|\mathcal{E}|}y_k^*tr(E_k\Phi_1^*),\quad t^*=\sum_{k=1}^{|\mathcal{E}|}y_k^*tr(E_k\Phi_2^*),
		\end{aligned}
		\label{eq:KKT_cond1}
	\end{equation}

	If there exists an optimal edge-weight realization of $\mathcal{\hat{G}(V,\hat{E},\hat{W})}$, $w_{|\mathcal{E}|+1}^*<0$, then the following optimization problem should have the same optimal value $\hat{t}^*$ as that of (\ref{eq:cvx_pro}), i.e., $t^*=\hat{t}^*$:
	\begin{equation}
		\begin{aligned}
			& \underset{\hat{t},\hat{y}_k,\hat{y}_0}{\text{minimize}}
			& & \hat{t}\\
			& \text{subject to}
			& & I-\hat{y}_0\bm{1}\bm{1}^T\preceq\sum_{k=1}^{|\mathcal{E}|+1}\hat{y}_kE_k\preceq \hat{t}I,\\
			& & & \hat{y}_{|\mathcal{E}|+1}\leq0.
		\end{aligned}
		\tag{P3}\label{eq:negative_weight_pro}
	\end{equation}
	Slater's condition also holds also for this problem. We will verify that there exists a solution where $\hat{y}_{|\mathcal{E}|+1}^*<0$ such that the KKT conditions for problem (\ref{eq:negative_weight_pro}) is satisfied while $t^*=\hat{t}^*$.

	The KKT conditions for Problem (\ref{eq:negative_weight_pro}) are
	\begin{equation}
		\begin{aligned}
			&tr(\hat{\Phi}_2^*)=1,\quad tr(\bm{1}\bm{1}^T\hat{\Phi}_1^*)=0,\\
			&tr(E_k\hat{\Phi}_2^*)-tr(E_k\hat{\Phi}_1^*)=0,\:k=1,\cdots,|\mathcal{E}|,\\
			&tr(E_{|\mathcal{E}|+1}\hat{\Phi}_2^*)-tr(E_{|\mathcal{E}|+1}\hat{\Phi}_1^*)\leq0,\\
			&\hat{y}_{|\mathcal{E}|+1}^*(tr(E_{|\mathcal{E}|+1}\hat{\Phi}_2^*)-tr(E_{|\mathcal{E}|+1}\hat{\Phi}_1^*))=0\\
			&tr(\hat{\Phi}_1^*)=\sum_{k=1}^{|\mathcal{E}|}\hat{y}_k^*tr(E_k\hat{\Phi}_1^*)+\hat{y}_{|\mathcal{E}|+1}^*tr(E_{|\mathcal{E}|+1}\hat{\Phi}_1^*),\\
			&t^*=\hat{t}^*=\sum_{k=1}^{|\mathcal{E}|}\hat{y}_k^*tr(E_k\hat{\Phi}_2^*)+\hat{y}_{|\mathcal{E}|+1}^*tr(E_{|\mathcal{E}|+1}\hat{\Phi}_2^*).
		\end{aligned}
		\label{eq:KKT_cond2}
	\end{equation}

	Compare (\ref{eq:KKT_cond1}) with (\ref{eq:KKT_cond2}) and consider the case $\Phi_1^*=\hat{\Phi}_1^*,\Phi_2^*=\hat{\Phi}_2^*$. The KKT conditions (\ref{eq:KKT_cond2}) are trivially fulfilled except potentially for the last two. Since strong duality holds for problem (\ref{eq:negative_weight_pro}), $tr(\hat{\Phi}_1^*)=\hat{t}^*$, and by $tr(E_k\hat{\Phi}_2^*)=tr(E_k\hat{\Phi}_1^*)$, the last two conditions are equivalent in this case. It is thus sufficient to consider only one of them.
	Now note that since $t^*\geq 1$ and $t^*=\sum_{k=1}^{|\mathcal{E}|}y_k^*tr(E_k\Phi_2^*)$, then not all $tr(E_k\Phi_2^*)$ equal zero. Therefore, if one chooses $\hat{y}_{|\mathcal{E}|+1}^*<0$, it is always possible to find a set of solution $\{\hat{y}_k^*\},\:k=1,\cdots,|\mathcal{E}|$ such that $t^*=\hat{t}^*=\sum_{k=1}^{|\mathcal{E}|}\hat{y}_k^*tr(E_k\hat{\Phi}_2^*)+\hat{y}_{|\mathcal{E}|+1}^*tr(E_{|\mathcal{E}|+1}\hat{\Phi}_2^*)$ holds.
\qed
\end{pf}

\begin{exmp}
	Consider the graph in Fig. \ref{fig:graph_ex}.
	\begin{figure}[!htpb]
		\centering
		\includegraphics[width=0.3\textwidth]{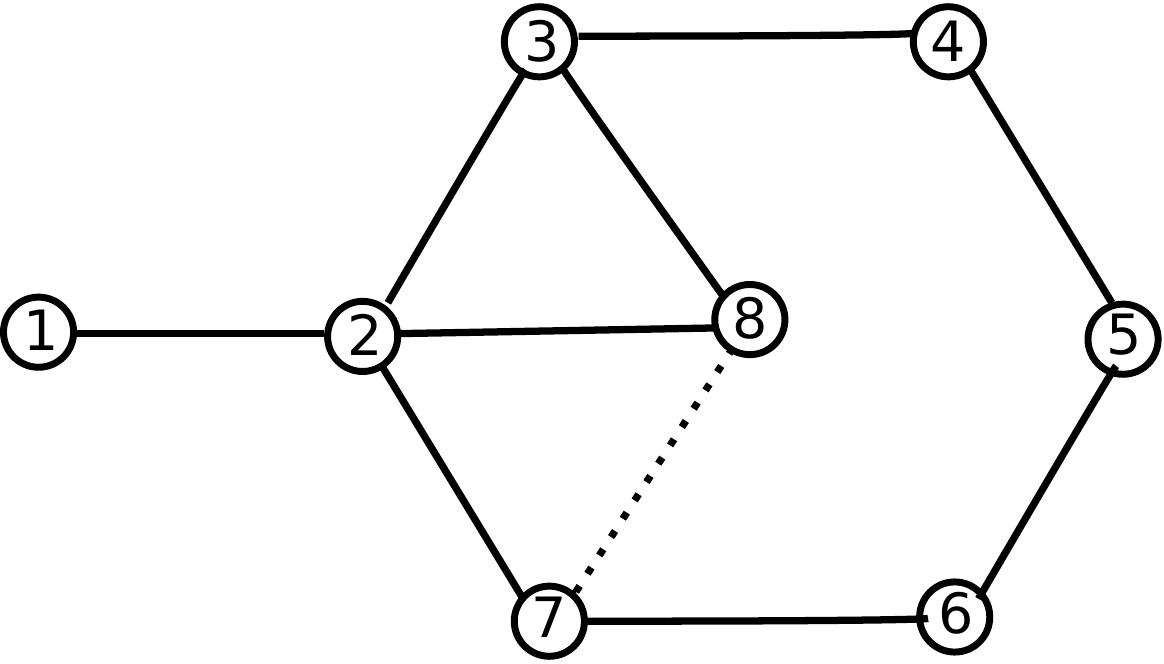}
		\caption{Example 1}
		\label{fig:graph_ex}
	\end{figure}
\end{exmp}
	We call the graph without the dashed line $\mathcal{G(V,W,E)}$ and the one with dashed line $\mathcal{\hat{G}}(V,\hat{W},\hat{E})$. One optimal edge-weight realization for $\mathcal{G(V,W,E)}$ is $y_{12}=2.0620$, $y_{23}=0.7215$, $y_{27}=1.4877$, $y_{28}=0.7667$, $y_{34}=2.0620$, $y_{38}=1.7297$, $y_{45}=1.3747$, $y_{56}=1.3747$, $y_{67}=2.0620$ with the optimal value $t^*=7.2480$. 

	One optimal edge-weight realization for $\mathcal{\hat{G}(V,\hat{W},\hat{E})}$ is $\hat{y}_{12}=2.0620$, $\hat{y}_{23}=1.0825$, $\hat{y}_{27}=1.0827$, $\hat{y}_{28}=0.8457$, $\hat{y}_{34}=2.0620$, $\hat{y}_{38}=0.9147$, $\hat{y}_{45}=1.3747$, $\hat{y}_{56}=1.3747$, $\hat{y}_{67}=2.0620,$ $\hat{y}_{78}= 0.9143$ with all edge-weight positive and same optimal value $\hat{t}^*=t^*=7.2480$ as that of $\mathcal{G(V,E,W)}$. 
	
	Another optimal edge-weight realization of $\mathcal{\hat{G}(V,\hat{W},\hat{E})}$ is $\hat{y}_{12}=2.0620$, $\hat{y}_{23}=0.6602$, $\hat{y}_{27}=1.5079$, $\hat{y}_{28}=0.8404$, $\hat{y}_{34}=2.0620$, $\hat{y}_{38}=1.8680$, $\hat{y}_{45}=1.3747$, $\hat{y}_{56}=1.3747$, $\hat{y}_{67}=2.0620$, $\hat{y}_{78}=-0.0456$. In this optimal edge-weight realization, agent 7 and 8 are ``competing with'' or ``push away from'' each other, but this will not effect the system performance.
	\qed

Another sufficient condition for the existence of negative edge weights is the following:
\begin{thm}\label{thm:negative weight}
	Suppose $\Phi_1^*$ and $\Phi_2^*$ are the optimal solutions to the dual problem (\ref{eq:dual_pro_simp}) for the graph $\mathcal{G(V,E,W)}$. If $tr(E_{|\mathcal{E}|+1}\Phi_2^*)-tr(E_{|\mathcal{E}|+1}\Phi_1^*)>0$, then the optimal solution to the optimal weight realization of $\mathcal{\hat{G}(V,\hat{E},\hat{W})}$ must satisfy $w_{|\mathcal{E}|+1}^*<0$ and $t^*>\hat{t}^*$.
\end{thm}
\begin{pf}
	To prove this, we first construct the following optimization problem:
	\begin{equation}
		\begin{aligned}
			& \underset{\hat{t},\hat{y}_k,\hat{y}_0}{\text{minimize}}
			& & \hat{t}\\
			& \text{subject to}
			& & I-\hat{y}_0\bm{1}\bm{1}^T\preceq\sum_{k=1}^{|\mathcal{E}|+1}\hat{y}_kE_k\preceq \hat{t}I,\\
			& & & \hat{y}_{|\mathcal{E}|+1}\geq0.
		\end{aligned}
		\tag{P4}\label{eq:positive_weight_pro}
	\end{equation}
	whose KKT-condition looks like:
	\begin{equation}
		\begin{aligned}
			&tr(\hat{\Phi}_2^*)=1,\quad tr(\bm{1}\bm{1}^T\hat{\Phi}_1^*)=0,\\
			&tr(E_k\hat{\Phi}_2^*)-tr(E_k\hat{\Phi}_1^*)=0,\:k=1,\cdots,|\mathcal{E}|,\\
			&tr(E_{|\mathcal{E}|+1}\hat{\Phi}_2^*)-tr(E_{|\mathcal{E}|+1}\hat{\Phi}_1^*)\geq0,\\
			&\hat{y}_{|\mathcal{E}|+1}^*(tr(E_{|\mathcal{E}|+1}\hat{\Phi}_2^*)-tr(E_{|\mathcal{E}|+1}\hat{\Phi}_1^*))=0,\\
			&tr(\hat{\Phi}_1^*)=\sum_{k=1}^{|\mathcal{E}|}\hat{y}_k^*tr(E_k\hat{\Phi}_1^*)+\hat{y}_{|\mathcal{E}|+1}^*tr(E_{|\mathcal{E}|+1}\hat{\Phi}_1^*),\\
			&t^*=\hat{t}^*=\sum_{k=1}^{|\mathcal{E}|}\hat{y}_k^*tr(E_k\hat{\Phi}_2^*)+\hat{y}_{|\mathcal{E}|+1}^*tr(E_{|\mathcal{E}|+1}\hat{\Phi}_2^*).
		\end{aligned}
		\label{eq:KKT_cond3}
	\end{equation}
	The optimal weight realization to $\mathcal{G(V,E,W)}$ can be written as $(\{y_k^*\},y_{|\mathcal{E}|+1}^*=0)$, where $k=1,\cdots,|\mathcal{E}|$, since not having the optimization variable $y_{|\mathcal{E}|+1}$ in (\ref{eq:cvx_pro}) is equivalent to $y_{|\mathcal{E}|+1}=0$). 
	Therefore, if $tr(E_{|\mathcal{E}|+1}\Phi_2^*)-tr(E_{|\mathcal{E}|+1}\Phi_1^*)>0$ holds, then it means that the optimal edge realization $(\{y_k^*\},y_{|\mathcal{E}|+1}^*=0)$ for $\mathcal{G(V,E,W)}$ is also optimal to (\ref{eq:positive_weight_pro}).

	But $(\{y_k^*\},y_{|\mathcal{E}|+1}^*=0)$ is not an optimal weight realization to $\mathcal{\hat{G}(V,\hat{E},\hat{W})}$ since it must satisfy the KKT conditions:
	\begin{align}
		tr(E_{|\mathcal{E}|+1}\hat{\Phi}_2^*)-tr(E_{|\mathcal{E}|+1}\hat{\Phi}_1^*)=0,
		\label{eq:KKT_cond4}
	\end{align}
	while by assumption it holds that $tr(E_{|\mathcal{E}|+1}\Phi_2^*)-tr(E_{|\mathcal{E}|+1}\Phi_1^*)>0$.
	Since $(\{y_k^*\},y_{|\mathcal{E}|+1}^*=0)$ is the optimal solution to (\ref{eq:positive_weight_pro}), this means $(\{y_k^*\},y_{|\mathcal{E}|+1}^*=0)$ is the best edge weight realization among all $(\{y_k\},y_{|\mathcal{E}|+1}\geq 0)$. But it is not the optimal edge weight realization to $\mathcal{\hat{G}(V,\hat{E},\hat{W})}$, which means for $\mathcal{\hat{G}(V,\hat{E},\hat{W})}$, it must satisfy $w_{|\mathcal{E}|+1}^*<0$ and $t^*>\hat{t}^*$.\qed
\end{pf}

\begin{exmp}
	Consider the graph in Fig. \ref{fig:graph_ex1}.
	\begin{figure}[!htpb]
		\centering
		\includegraphics[width=0.3\textwidth]{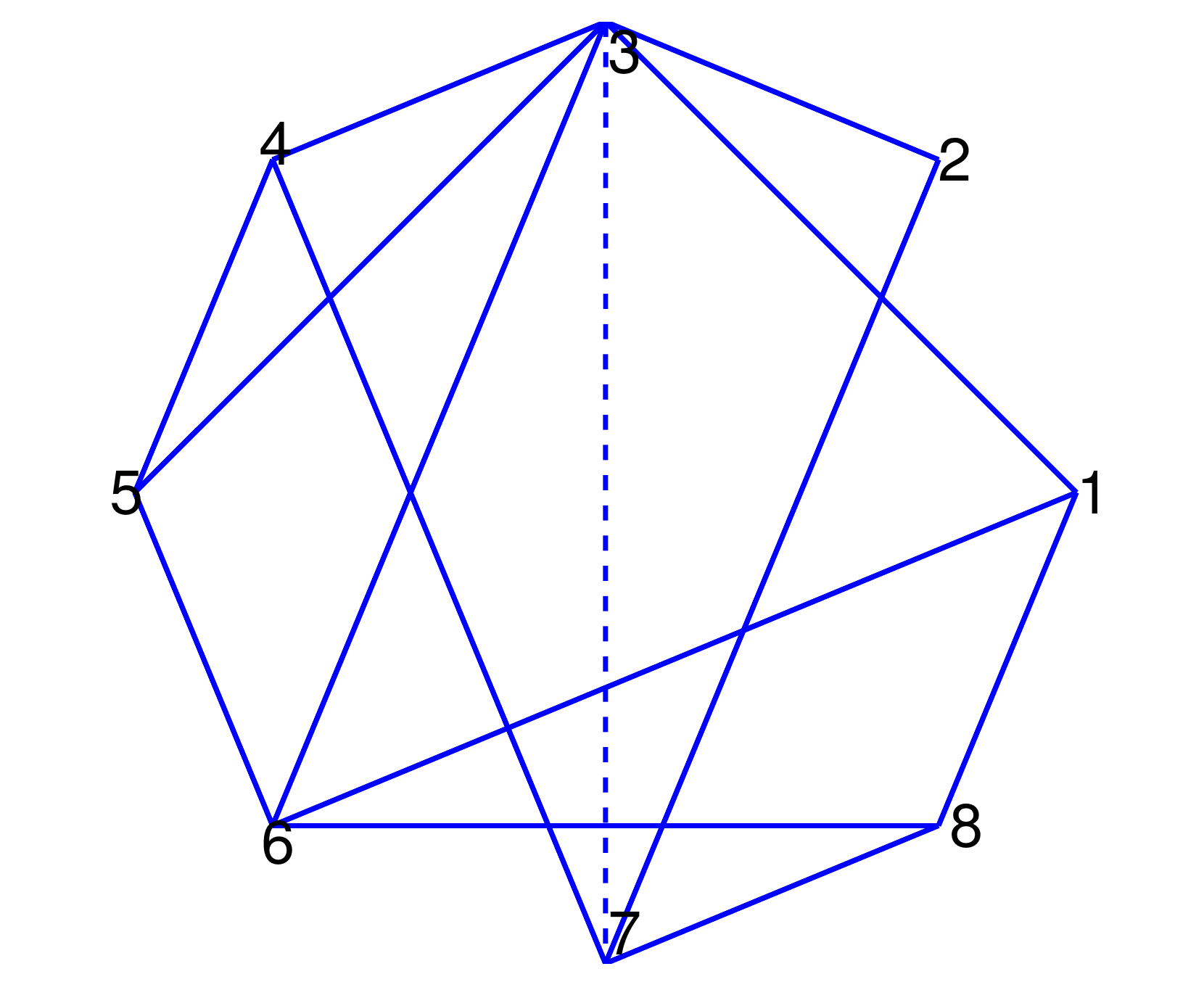}
		\caption{Example 2}
		\label{fig:graph_ex1}
	\end{figure}
\end{exmp}
	We call the graph without the dashed line $\mathcal{G(V,W,E)}$ and the one with dashed line $\mathcal{\hat{G}(V,\hat{W},\hat{E})}$. One optimal edge-weight realization for $\mathcal{G(V,W,E)}$ is $y_{13}=0.5329$, $y_{16}=0.4896$, $y_{18}=0.5675$, $y_{23}=0.7777$, $y_{27}=0.5072$, $y_{34}=0.1534$, $y_{35}=0.4235$, $y_{36}=0.4331$, $y_{45}=0.7028$, $y_{47}=0.5991$, $y_{56}=0.7604$, $y_{68}=0.4663$, $y_{68}=0.8208$ with the optimal value $t^*=3.0592$.  $tr(E_{37}\Phi_2^*)-tr(E_{37}\Phi_1^*)=0.0467$.
	
	One optimal edge-weight realization for $\mathcal{\hat{G}(V,\hat{W},\hat{E})}$ is $\hat{y}_{13}=0.5317$, $\hat{y}_{16}=0.4789$, $\hat{y}_{18}=0.6191$, $\hat{y}_{23}=0.7794$, $\hat{y}_{27}=0.4981$, $\hat{y}_{34}=0.1581$, $\hat{y}_{35}=0.4317$, $\hat{y}_{36}=0.4368$, $\hat{y}_{37}=-0.0495$, $\hat{y}_{45}=0.7112$, $\hat{y}_{47}=0.5932$, $\hat{y}_{56}=0.7577$, $\hat{y}_{68}=0.4658$, $\hat{y}_{68}=0.8183$ with the optimal value $\hat{t}^*=3.0581$. 
	
	Moreover, to illustrate negative weights may lead to a lower control energy cost, choose the matrices in \eqref{eq:system} and ARE as
\eqref{eq:A B Q} and $\varepsilon=0$.
The initial value $\hat{X}_0$ (recall (LQ2) in the manuscript) is chosen as $\hat{X}_0=[ 1.4090, 1.4172, 0.6715, -1.2075,$ $0.7172, 1.6302, 0.4889, 1.0347, 0.7269, -0.3034, 0.2939,$ $-0.7873, 0.8884, -1.1471]^T$. The actual control energy cost under the allowance of negative edge weight is $156.3912$, while the actual control energy cost under the constraint of non-negative edge weight is $156.4276$. (The same initial value is used for both cases.)
	\qed

\begin{rem}
From the two theorems above, it can be seen that if we already have an optimal solution of the original graph and we want to add an edge between two nodes, we can use the sign of $tr(E_{|\mathcal{E}|+1}\Phi_2^*)-tr(E_{|\mathcal{E}|+1}\Phi_1^*)$ to determine if the two agents will "compete with" each other or not before doing the computation. Moreover, if $tr(E_{|\mathcal{E}|+1}\Phi_2^*)-tr(E_{|\mathcal{E}|+1}\Phi_1^*)<0$, then the optimal solution to the optimal weight realization of $\mathcal{\hat{G}(V,\hat{E},\hat{W})}$ must satisfy $w_{|\mathcal{E}|+1}^*>0$ and $t^*>\hat{t}^*$, which can be shown in a similar way as Theorem \ref{thm:negative weight}.
\end{rem}

\subsection{Distributed Optimization in the Case of Regular Graphs}\label{sec:distri_opt}
When the graph is $\kappa$-regular, namely, each node has $\kappa$ neighbours, this also means $|\mathcal{E}(i)|=\kappa,\forall i$; we are able to solve the problem (\ref{eq:cvx_pro}) in a distributed manner. In fact, this case is quite common in the real applications since usually for homogeneous nodes that construct the network, each node has the same amount of communication channel. This assumption will be used later and explained in Remark \ref{rem:kappa regular}.
\subsubsection{Intersection of Convex Sets}
 In order to develop the distributed optimization algorithm for regular graphs, we first rewrite the dual problem (\ref{eq:dual_pro_simp}) into a standard SDP form:
\begin{equation}
\begin{aligned}
		& \underset{\Phi}{\text{maximize}}
		& & tr(\hat{I}\Phi)\\
		& \text{subject to}
		& & tr(\tilde{I}\Phi)=1,\\
		& & & tr(\mathbbm{1}\Phi)=0,\\
		& & & tr(\hat{E}_k\Phi)=0,\: k=1,2,\cdots,|\mathcal{E}|,\\
		& & & \Phi\succeq 0.\\
	\end{aligned}
	\tag{DS}\label{eq:standard dual}
\end{equation}
where $\Phi=diag(\Phi_1,\Phi_2)$, $\hat{I}=diag(I_N,\bm{0})$, $\tilde{I}=diag(\bm{0},I_N)$, $\mathbbm{1}=diag(\bm{1}\bm{1}^T,\bm{0})$ and $\hat{E}_k=diag(-E_k,E_k)$.

From the node's perspective, we can rewrite the feasible domain $\Omega$ of (\ref{eq:standard dual}) as an intersection of node's local feasible domain:
\begin{equation}\nonumber
\begin{aligned}
\Omega&=\big\{\Phi\:|\:tr(\tilde{I}\Phi)=1,\:tr(\mathbbm{1}\Phi)=0,\:\Phi\succeq 0,\\&tr(\hat{E}_k\Phi)=0,\: k=1,2,\cdots,|\mathcal{E}|,\big\}\\
&=\bigcap_{i\in\mathcal{V}}\Omega_i,
\end{aligned}
\end{equation}
where
\begin{equation}
\begin{aligned}
\Omega_i=\big\{\Phi^{(i)}\:|\:tr(\tilde{I}\Phi^{(i)})=1,\:tr(\mathbbm{1}\Phi^{(i)})=0,\\\Phi^{(i)}\succeq 0,\:tr(\hat{E}_k\Phi^{(i)})=0,\: k\in\mathcal{E}(i)\big\}.
\end{aligned}
\end{equation}
Therefore the idea of this up-coming distributed optimization set-up is that each node solves its own $\Phi^{(i)*}$ and all $\Phi^{(i)*}$ reach consensus in the intersection of the $\Omega_i$. Hence a penalty term on consensus of $\Phi^{(i)}$ is added to the objective function and leads to the following problem
\begin{equation}
\begin{aligned}
		& \underset{\Phi^{(i)}}{\text{minimize}}
		& & \sum_{i\in\mathcal{V}}\big\{-\frac{1}{N}tr(\hat{I}\Phi^{(i)})\\&&&+\frac{M}{2}\sum_{j\in\mathcal{N}(i)}\|\Phi^{(i)}-\Phi^{(j)}\|_F^2\big\}\\
		& \text{subject to}
		& & \Phi^{(i)}\in \Omega_i,\\
	\end{aligned}
	\label{eq:distributed dual1}
\end{equation}
where $M$ is the penalty parameter.
\begin{prop}
Problem (\ref{eq:distributed dual1}) has the following solution $\big\{\Phi^{(1)*},\cdots,\Phi^{(N)*}\big\}$, where $\Phi^{(i)*}\approx\Phi^{(j)*},\:\forall i\neq j$ and $\Phi^{(i)*}$ is almost the optimal solution to (\ref{eq:standard dual}).
\end{prop}
\begin{pf}
The KKT conditions of (\ref{eq:distributed dual1}) read
\begin{align}
&tr(\tilde{I}\Phi^{(i)*})=1,\:tr(\mathbbm{1}\Phi^{(i)*})=0,\: \Phi^{(i)*}\succeq 0,\nonumber\\
&tr(\hat{E}_k\Phi^{(i)*})=0,\: k\in\mathcal{E}(i),\:S^{(i)*}\succeq 0,\nonumber\\
&-\frac{1}{N}\hat{I}+t^{(i)*}\tilde{I}+y_0^{(i)*}\mathbbm{1}-\sum_{k\in\mathcal{E}(i)}y_k^{(i)*}\hat{E}_k-S^{(i)*}\nonumber\\
&+2M\sum_{j\in\mathcal{N}(i)}(\Phi^{(i)*}-\Phi^{(j)*})=0,\label{eq:distributed KKT1}\\
&tr(S^{(i)*}\Phi^{(i)*})=0, \forall i\in \mathcal{V},\nonumber
\end{align}
Because of the penalty term in (\ref{eq:distributed dual1}), $\Phi^{(i)*}\approx\Phi^{(j)*},\:\forall i\neq j$, then summing-up the second last equations in (\ref{eq:distributed KKT1}) for all $i\in\mathcal{V}$ will result in "almost" the same KKT conditions of (\ref{eq:standard dual}).
\end{pf}
\begin{rem}
We would like to use Primal Dual Interior Point Method (PDIPM) (see \cite{alizadeh1998primal}) to solve the SDP problem because it is a robust and efficient method for SDP problems. The reason we use a penalty term on consensus of $\Phi^{(i)}$ instead of using consensus constraint $\sum_{j\in\mathcal{N}(i)}(\Phi^{(i)}-\Phi^{(j)})=0$ is that PDIPM uses Newton iteration to solve the perturbed KKT conditions. However, if the consensus constraint is introduced, then this also results in the introduction of Lagrangian multiplier terms $\sum_{j\in\mathcal{N}(i)}(Z^{(i)}-Z^{(j)})$ into the perturbed KKT conditions. Together with the consensus contraint $\sum_{j\in\mathcal{N}(i)}(\Phi^{(i)}-\Phi^{(j)})=0$, a rank deficient matrix $L_u\otimes I_N$ will show up in the matrix $\mathcal{Q}^l$ while equation $\mathcal{Q}^l\Delta p^l=r^l$ solves the Newton direction $\Delta p^l$ of the $l$th iteration and hence makes $\mathcal{Q}^l$ singular.\qed
\end{rem}
\subsubsection{Primal Dual Interior Point Method}\label{sec:primal_dual_interior}
We solve problem (\ref{eq:distributed dual1}) using PDIPM. If we are not satisfied with result on the consensus among $\Phi^{(i)^*}$, then $M$ is increased and we use the optimal solution as the initial guess for the new problem, which is called as "warm start". This technique is called Sequential Unconstrained Minimization Technique (SUMT).

Now we solve the problem (\ref{eq:distributed dual1}) using PDIPM. The penalized barrier problem of parameter $\rho>0$ reads
\begin{equation}\nonumber
\begin{aligned}
		& \underset{\Phi^{(i)}}{\text{minimize}}
		& & \sum_{i\in\mathcal{V}}\big\{-\frac{1}{N}tr(\hat{I}\Phi^{(i)})-\rho\ln \det \Phi^{(i)}\\
		&&&+\frac{M}{2}\sum_{j\in\mathcal{N}(i)}\|\Phi^{(i)}-\Phi^{(j)})\|_F^2\big\}\\
		& \text{subject to}
		& & \Phi^{(i)}\in \Omega_i.\\
	\end{aligned}
\end{equation}

In order to make the matrix $\mathcal{Q}^l$ symmetric in the $\mathcal{Q}^l\Delta p^l=r^l$ while computing the Newton step $\Delta p^l$, the perturbed complementary slackness $\Phi^{(i)}S^{(i)}=\rho I_{2N}$ is rewritten as $\Phi^{(i)}S^{(i)}+S^{(i)}\Phi^{(i)}=2\rho I_{2N}$ in the central path. Hence the Newton step of the $l$th iteration is computed through solving
\begin{equation}
\begin{aligned}
&tr(\tilde{I}\Delta\Phi^{(i)l})=1-tr(\tilde{I}\Phi^{(i)l}),\:tr(\mathbbm{1}\Delta\Phi^{(i)l})=-tr(\mathbbm{1}\Phi^{(i)l}),\\
&tr(\hat{E}_k\Delta\Phi^{(i)l})=-tr(\hat{E}_k\Phi^{(i)l}),\\
&2M\sum_{j\in\mathcal{N}(i)}(\Delta\Phi^{(i)l}-\Delta\Phi^{(j)l})+\Delta t^{(i)l}\tilde{I}+\Delta y_0^{(i)l}\mathbbm{1}\\
&-\sum_{k\in\mathcal{E}(i)}\Delta y_k^{(i)l}\hat{E}_k+\Delta S^{(i)l}=r_{primal}^{(i)l},\\
&\Delta\Phi^{(i)l}S^{(i)l}+\Phi^{(i)l}\Delta S^{(i)l}+\Delta S^{(i)l}\Phi^{(i)l}+S^{(i)l}\Delta\Phi^{(i)l}\\
&=\rho I_{2N}-\Phi^{(i)l}S^{(i)l}+ S^{(i)l}\Phi^{(i)l},
\end{aligned}
\label{eq:newton step}
\end{equation}
where $r_{primal}^{(i)l}=\frac{1}{N}-2M\sum_{j\in\mathcal{N}(i)}(\Phi^{(i)l}-\Phi^{(j)l})-t^{(i)l}\tilde{I}-y_0^{(i)l}\mathbbm{1}+\sum_{k\in\mathcal{E}(i)}y_k^{(i)l}\hat{E}_k-S^{(i)l}$. We do symmetric vectorization on both hand sides of (\ref{eq:newton step}). Note that when we do the symmetric vectorization on the last equation, it holds that 
\begin{equation}\nonumber
\begin{aligned}
&\Delta\Phi^{(i)l}S^{(i)l}+\Phi^{(i)l}\Delta S^{(i)l}+\Delta S^{(i)l}\Phi^{(i)l}+S^{(i)l}\Delta\Phi^{(i)l}\\
&=\rho I_{2N}-\Phi^{(i)l}S^{(i)l}+ S^{(i)l}\Phi^{(i)l}\\
\Leftrightarrow&(S^{(i)l}\otimes_sI_{2N})svec(\Delta\Phi^{(i)l})+(\Phi^{(i)l}\otimes_sI_{2N})svec(\Delta S^{(i)l})\\
&=svec(\rho I_{2N}-\frac{1}{2}(\Phi^{(i)l}S^{(i)l}+S^{(i)l}\Phi^{(i)l})).
\end{aligned}
\end{equation}
Therefore (\ref{eq:newton step}) of the $l$th iteration can be written in the matrix form:
\begin{equation}
 \mathcal{Q}^l\Delta p^l=r^l,
 \label{eq:newton_step_mat_form}
\end{equation}
where $\mathcal{Q}^l=diag(\mathcal{Q}^{(1)l},\cdots,\mathcal{Q}^{(N)l})+L_u\otimes \mathcal{C}=[\bar{\mathcal{Q}}^{(1)l^T},\cdots,\bar{\mathcal{Q}}^{(N)l^T}]^T$, $\Delta p^l=[\Delta p^{(1)l^T},\cdots,\Delta p^{(N)l^T}]^T$, $r^l=[r^{(1)l^T},\cdots,r^{(N)l^T}]^T$. (The structure of the detailed matrices can be found in the appendix.)
\begin{rem}\label{rem:kappa regular}
Although for each agent $i$, it only knows its own matrix block $\mathcal{Q}^{(i)l}$, it does not mean that agent $i$ is ignorant of the column dimension of the row block $\bar{\mathcal{Q}}^{(i)l}$. Since the graph is $\kappa$ regular, each block $\mathcal{Q}^{(i)l}$ must have the same dimension. Thus each row block $\bar{\mathcal{Q}}^{(i)l}$ is totally accessible by agent $i$.\qed 
\end{rem}

Gauss elimination can be further done on $svec(\Delta S^{(i)l})$ in (\ref{eq:newton_step_mat_form}) and hence we get the equation $\hat{\mathcal{Q}}^l\Delta\hat{p}^l=\hat{r}^l$, where $\hat{\mathcal{Q}}^l=diag(\hat{\mathcal{Q}}^{(1)l},\cdots,\hat{\mathcal{Q}}^{(N)l})+L_u\otimes \hat{\mathcal{C}}=[\tilde{\mathcal{Q}}^{(1)l^T},\cdots,\tilde{\mathcal{Q}}^{(N)l^T}]^T$, $\Delta\hat{p}^l=[\Delta\hat{p}^{(1)l^T},\cdots,\Delta\hat{p}^{(N)l^T}]^T$, $\hat{r}^l=[\hat{r}^{(1)lT},\cdots,\hat{r}^{(N)lT}]^T$. Now that every node knows the block $\tilde{\mathcal{Q}}^{(i)l}$ and $\hat{r}^{(i)l}$, the approach proposed by \cite{mou2015distributed} is used to solve the system of linear equations in a distributed manner. By using the algorithm, each node is able to compute its own $\Delta \hat{p}^l$ and all $\Delta \hat{p}^l$ computed by each node will reach consensus. It has been proved that the algorithm converge exponentially. The next step is to compute the step length $\alpha^l$. The step length $\alpha^{(i)l}$ for each node is computed through
\begin{equation}\nonumber
\begin{aligned}
\alpha^{(i)l}=&\max\big\{\alpha^{(i)l}\:|\:\Phi^{(i)l}+\alpha^{(i)l}\Delta\Phi^{(i)l}\succeq 0,\\
&S^{(i)l}+\alpha^{(i)l}\Delta S^{(i)l}\succeq 0\big\}.
\end{aligned}
\end{equation}
Now every node has its own $\alpha^{(i)l}$, but in order to make a correct Newton step, we need an $\alpha^l$ such that  $\Phi^{(i)l}+\alpha^l\Delta\Phi^{(i)l}\succeq 0$, $S^{(i)l}+\alpha^l\Delta S^{(i)l}\succeq 0$ for all nodes. In order to do this, each node can compare its own $\alpha^{(i)l}$ with those of its neighbours; and choose
\begin{equation}\nonumber
\alpha^{(i)l}\leftarrow \min\big\{\alpha^{(i)l},\alpha^{(j)l}\big\},\quad \forall j\in \mathcal{N}(i).
\end{equation}
This process can be seen as the propagation of the smallest step-length $\alpha^{*l}$ from its node to the entire network. The maximum number of steps for every node to get $\alpha^{*l}$ is the length of the longest path in the graph.

As a summary, for each node the algorithm can be expressed as:
\begin{alg}
Distributed PDIPM\\
\begin{tabular}{l}\hline
\textbf{Input}: $\Phi^{(i)0}$, $\forall i\in \mathcal{V}$, $\epsilon>0$, $\theta\in(0,1)$, tol, $\xi>1$\\
\hline
\textbf{While} $\max\|\Phi^{(i)*}-\Phi^{(j)*}\|_F>$ tol, $j\in\mathcal{N}(i)$\\
\hspace{4ex}\textbf{While} $\rho>\epsilon$\\
\hspace{8ex}each node compute $\Delta \hat{p}^l$\\
\hspace{8ex}compute $\alpha^{*l}$\\
\hspace{8ex}$\Phi^{(i)(l+1)}\leftarrow\Phi^{(i)l}+\alpha^{*l}\Delta\Phi^{(i)l}$\\
\hspace{8ex}$S^{(i)(l+1)}\leftarrow S^{(i)l}+\alpha^{*l}\Delta S^{(i)l}$\\
\hspace{8ex}$\rho\leftarrow\rho\theta$\\
\hspace{8ex}$l\leftarrow l+1$\\
\hspace{4ex}\textbf{end}\\
\hspace{4ex}$\Phi^{(i)*}\leftarrow \Phi^{(i)l}$\\
\hspace{4ex}$M\leftarrow M\xi$\\
\textbf{end}
\end{tabular}
\end{alg}
\begin{exmp}
Consider the graph in (\ref{fig:ring}).
\begin{figure}[!htpb]
		\centering
		\includegraphics[width=0.25\textwidth]{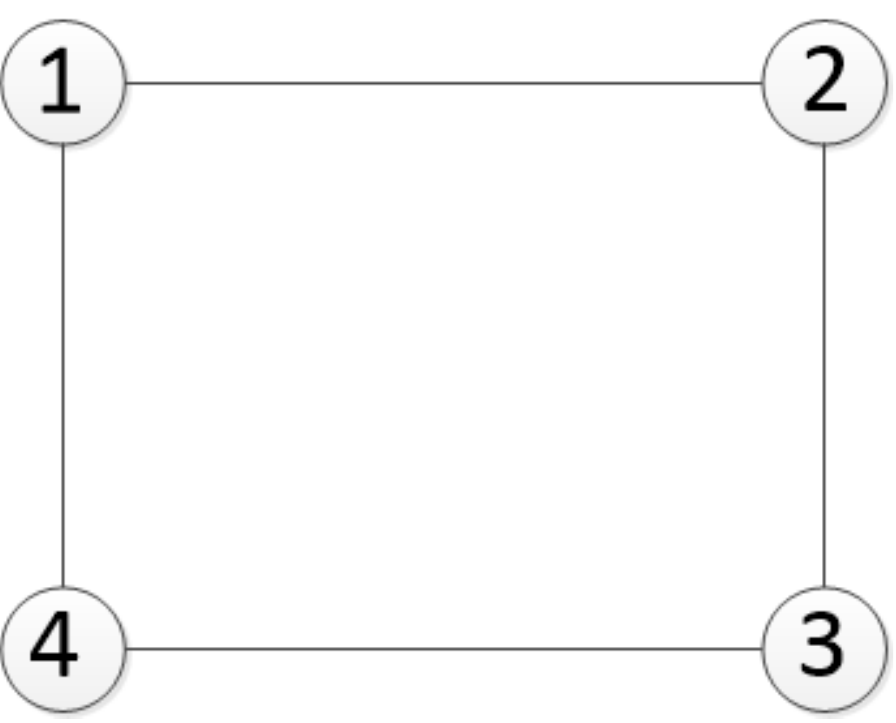}
		\caption{Example 3}
		\label{fig:ring}
\end{figure}
Choose $M=500$, $\epsilon=10^{-13}$, $\theta=0.1$, $\rho=10000$, $\xi=2$ and we get $y_{12}^{(1)*}=y_{12}^{(2)*}=y_{23}^{(2)*}=y_{23}^{(3)*}=y_{34}^{(3)*}=y_{34}^{(4)*}=y_{41}^{(4)*}=y_{41}^{(1)*}=0.2505$, $\max\|\Phi^{(i)*}-\Phi^{(j)*}\|_F=1.0016\times 10^{-4},j\in\mathcal{N}(i),\forall i\in \mathcal{V}$. Hence $y_{12}^*=y_{23}^*=y_{34}^*=y_{41}^*=0.5010$ while the optimal solution of centralized problem (\ref{eq:cvx_pro}) is $y_{12}^*=y_{23}^*=y_{34}^*=y_{41}^*=0.5000$.
\end{exmp}

\section{Conclusion and Future Work}\label{sec:conclusion}
In this paper, we study consensus control for linear systems with optimal energy cost. Due to the hardness of solving the problem, we focus on a classical family of controllers that are based on ARE and guarantee consensus. A suboptimal energy controller that depends only on the relative information between the agents has been constructed by optimizing the control gain as well as the edge weights of the graph. We have shown that the controller coincides with the results in \cite{thunberg2016optimal} when the graph topology is complete and provide two sufficient conditions for the existence of negative optimal edge weights when a new edge is added to the graph. Moreover, we proposed a distributed optimization algorithm for solving the problem when the graph is $\kappa$-regular.

The authors are currently working on improving the convergence rate of the distributed optimization algorithm as well as developing a distributed optimization algorithm for maximizing the algebraic connectivity of graphs by adjusting the edge-weights.
\bibliographystyle{model5-names}
\bibliography{sync_opt} 

\appendix
\section{The Structure of the Matrices in Newton Iteration}
\begin{equation}\nonumber
\begin{aligned}
 &\Delta p^{(i)l^T}=[\Delta y_0^{(i)l},t^{(i)l},\Delta y_{k_1}^{(i)l},\cdots,\Delta y_{k_\kappa}^{(i)l},svec(\Delta\Phi^{(i)l})^T,\\
 &svec(\Delta S^{(i)l})^T]=[\Delta p^{(i)l^T}_{primal},\Delta p^{(i)l^T}_{dual},\Delta p^{(i)l^T}_{comp}],
 \end{aligned}
\end{equation}
where $k_1,\cdots,k_\kappa\in \mathcal{E}(i)$.
\begin{equation}\nonumber
\begin{aligned}
\mathcal{Q}^{(i)}=\begin{bmatrix}
\bm{0}_{\kappa+2} &\mathcal{A}_i &\bm{0}_{(\kappa+2)\times N(2N+1)}\\
\mathcal{A}_i^T &\bm{0}_{N(2N+1)}& -I_{N(2N+1)} \\
\bm{0}_{(\kappa+2)\times N(2N+1)}^T &S^{(i)l}\otimes_s I_{2N} & \Phi^{(i)l}\otimes_s I_{2N} 
\end{bmatrix},
\end{aligned}
\end{equation}
\begin{equation}\nonumber
\begin{aligned}
\mathcal{C}=\begin{bmatrix}
\bm{0}_{\kappa+2} &\bm{0}_{(\kappa+2)\times N(2N+1)} &\bm{0}_{(\kappa+2)\times N(2N+1)}\\
\bm{0}_{(\kappa+2)\times N(2N+1)}^T &2MI_{N(2N+1)} & \bm{0}_{N(2N+1)} \\
\bm{0}_{(\kappa+2)\times N(2N+1)}^T &\bm{0}_{N(2N+1)}^T & \bm{0}_{N(2N+1)} 
\end{bmatrix},
\end{aligned}
\end{equation}
\begin{equation}\nonumber
\begin{aligned}
\mathcal{A}_i^T=[svec(\mathbbm{1}),svec(\tilde{I}),svec(\hat{E}_{k_1}),\cdots,svec(\hat{E}_{k_\kappa})],
\end{aligned}
\end{equation}
where $k_1,\cdots,k_\kappa\in \mathcal{E}(i)$.

\begin{equation}\nonumber
\begin{aligned}
&r^{(i)l^T}=[-tr(\mathbbm{1}\Phi^{(i)l}),1-tr(\tilde{I}\Phi^{(i)l}),tr(\hat{E}_{k_1}\Phi^{(i)l}),\cdots,\\
&tr(\hat{E}_{k_\kappa},\Phi^{(i)l}),r_{primal}^{(i)l^T},svec(\rho I_{2N}-\frac{1}{2}(\Phi^{(i)l}S^{(i)l}+S^{(i)l}\Phi^{(i)l}))^T]\\
&=[r_{dual}^{(i)l^T},r_{primal}^{(i)l^T},r_{comp}^{(i)l^T}],
\end{aligned}
\end{equation}
where $k_1,\cdots,k_\kappa\in \mathcal{E}(i)$.

\begin{equation}\nonumber
\hat{\mathcal{Q}}^{(i)l}=
\begin{bmatrix}
\bm{0}_{\kappa+2} &\mathcal{A}_i \\
\mathcal{A}_i^T &(\Phi^{(i)l}\otimes_s I_N)^{-1}(S^{(i)l}\otimes_s I_N)
\end{bmatrix},
\end{equation}
\begin{equation}\nonumber
\hat{\mathcal{C}}=
\begin{bmatrix}
\bm{0}_{\kappa+2} &\bm{0}_{(\kappa+2)\times N(2N+1)}\\
\bm{0}_{(\kappa+2)\times N(2N+1)}^T &2MI_{N(2N+1)}
\end{bmatrix},
\end{equation}
\begin{equation}\nonumber
\begin{aligned}
&\hat{p}^{(i)l^T}=[\Delta p^{(i)l^T}_{primal},\Delta p^{(i)l^T}_{dual}],\\
&\hat{r}^{(i)l^T}=[r_{dual}^{(i)l^T},r_{primal}^{(i)l^T}+r_{comp}^{(i)l^T}(\Phi^{(i)l}\otimes_s I_N)^{-T}].
\end{aligned}
\end{equation}
\end{document}